\documentclass[12pt,reqno,a4paper]{amsart} 

\usepackage{mathtools}
\usepackage[breaklinks]{hyperref} 

\usepackage[headings]{fullpage} 

\usepackage{pslatex} 

\usepackage{microtype}  
 
\renewcommand{\qedsymbol}{$\square$}
\newenvironment{Proof}[1][Proof]{\par\noindent\textbf{#1.}~}
{\hfill\qedsymbol\smallskip\par}

\newcommand{\dx}{\mathrm{d}}
\newcommand{\eps}{\varepsilon}

\newcommand{\Stilde}{\widetilde{S}}
\newcommand{\Etilde}{\widetilde{E}}

\newcommand{\Ttilde}{\widetilde{T}}
\newcommand{\second}{\prime\prime} 
 
\newcommand{\Odip}[2]{\mathcal{O}_{#1}\!\left(#2\right)\mathchoice{\!}{}{}{}}
\newcommand{\Odig}[1]{\mathcal{O}\Bigl(#1\Bigr)\mathchoice{\!}{}{}{}}

\newcommand{\Odim}[1]{\mathcal{O}\bigl(#1\bigr)}
\newcommand{\Odi}[1]{\Odip{}{#1}}
\newcommand{\odip}[2]{{o}_{#1}\!\left(#2\right)\mathchoice{\!}{}{}{}}
\newcommand{\odi}[1]{\odip{}{#1}}
\newcommand{\odim}[1]{o\bigl(#1\bigr)}
 
\makeatletter
\newenvironment{Biggcases}{%
  \matrix@check\Biggcases\env@Biggcases
}{%
  \endarray %
}
\def\env@Biggcases{%
  \let\@ifnextchar\new@ifnextchar
  \Biggl\lbrace
  \def\arraystretch{1.2}%
  \array{@{}l@{\quad}l@{}}%
}
\makeatother

\newtheoremstyle{sltheorems}
{10pt}
{6pt}
{\slshape}
{}
{\bfseries}
{.}
{.5em}
{\thmname{#1}\thmnumber{ #2}\thmnote{ (#3)}}

\theoremstyle{sltheorems} 
\newtheorem{Theorem}{Theorem}
\newtheorem{Lemma}{Lemma} 
\newtheorem{Remark}{Remark}

\allowdisplaybreaks
\begin{document}
\title[asymptotic formulae for binary problems, II: density $1$]{Short intervals asymptotic formulae for binary problems with primes and powers, II: density $1$} 
\author[]{Alessandro Languasco \lowercase{and} Alessandro Zaccagnini}

\begin{abstract}
We prove that suitable asymptotic formulae in short intervals hold for the problems
of representing an integer as a sum of a prime square and a square, or a prime square.
Such results are obtained both assuming the Riemann Hypothesis and
in the unconditional case.
\end{abstract}

\maketitle

\allowdisplaybreaks
   
\section{Introduction}
In this second paper devoted to study asymptotic formulae in short intervals for additive problems with primes and squares, we focus our attention on density-one problems, \emph{i.e.}, on representing integers as sum of two squares. 
We considered the case of the sum of a prime and a square in our paper \cite{LanguascoZ2015c}.

We will consider two separate cases depending on the number of prime squares involved in the summations.
Let $\eps>0$, $N$ be a sufficiently large integer and  let further $H$ be an integer such that $N^{\eps}<H=\odi{N}$ as $N \to \infty$.
Taking $n\in [N, N+H]$, the key quantities  are
\[
r^{\second}_{2,2}(n) = 
\sum_{p_{1}^{2}+p_{2}^2=n}
\log p_{1} \log p_{2}
\quad
\textrm{and}
\quad
r^{\prime}_{2,2}(n) = 
\sum_{p^{2}+m^2=n}
\log p. 
\]
Since it is well known that the expected behaviour of such functions is erratic,  
to work in a more regular situation  we will study their average asymptotics  over a suitable short interval.

We have the following results which extend and improve the ones cited in the Introduction of the paper by Daniel \cite{Daniel2001}. 
We write $f=\infty(g)$ for $g=\odi{f}$.
\begin{Theorem}
\label{asymp-part3}  Assume the Riemann Hypothesis (RH) holds. Then
\[
\sum_{n = N+1}^{N + H} r^{\second}_{2,2}(n)
=
 \frac{\pi}{4}   H
+ \Odig{\frac{H^2}{N} +  H^{1/2}N^{1/4}(\log N)^{3/2}}
\]
as  $N\to \infty$ uniformly for $\infty(N^{1/2}(\log N)^{3})\le H \le  \odi{N}$. 
\end{Theorem}

\begin{Theorem}
\label{asymp-unc-part3}
Let $\eps>0$. Then there exists a constant $C=C(\eps)>0$ such that
\[
\sum_{n = N+1}^{N + H} r^{\second}_{2,2}(n)
=
 \frac{\pi}{4} H
 + 
\Odig{H \exp \Big( - C  \Big( \frac{\log N}{\log \log N} \Big)^{1/3} \Big)}
\]
as  $N\to \infty$  uniformly for   $N^{7/12+\eps}\le H \le N^{1-\eps}$. 
\end{Theorem}
 
We remark that Plaksin \cite{Plaksin1981} (see Lemma 11 there)
proves the case $H=N$ of Theorem \ref{asymp-unc-part3}
with a stronger error term of the form
$N \exp ( - C  (\log N)^{1/2} )$. Following its proof it
is clear that it can be further improved to
$N \exp ( - C   (\log N)^{3/5}(\log \log N)^{-1/5} )$.
The comparative weakness of our error term is due to the
use of the zero-density estimates for the Riemann zeta-function
(we need them to be able to get a short interval result).
A direct trial following the lines of  Lemma 11 of Plaksin \cite{Plaksin1981}
leads to weaker uniformity ranges:  
$H\gg N^{3/4}L^A$, for some $A>0$, assuming RH and  $H\gg N^{7/24+1/2+\eps}$ unconditionally.
Here $L = \log N$.

Concerning the sum of a prime square and a square, we have
\begin{Theorem}
\label{asymp-part4}  Assume the Riemann Hypothesis  holds. Then
\[
\sum_{n = N+1}^{N + H} r^{\prime}_{2,2}(n)
=
\frac{\pi}{4} H 
+
 \Odig{\frac{H^2}{N}+\frac{H \log \log N} {(\log N)^{1/2}}  }
 \]
as  $N\to \infty$ uniformly for $\infty(N^{1/2}(\log N)^{2}) \le H \le \odi{N}$. 
\end{Theorem}

\begin{Theorem}
\label{asymp-unc-part4}
Let $\eps>0$. Then there exists a constant $C=C(\eps)>0$ such that
\[
\sum_{n = N+1}^{N + H} \, r^{\prime}_{2,2}(n)
=
\frac{\pi}{4} H 
+ 
\Odig{  H \exp \Big( - C  \Big( \frac{\log N}{\log \log N} \Big)^{1/3} \Big)}
\]
as  $N\to \infty$ uniformly for   $N^{7/12+\eps}\le H \le N^{1-\eps}$. 
\end{Theorem}

An argument similar to the proof Lemma 11 of  Plaksin \cite{Plaksin1981}  
proves the case $H=N$ of Theorem \ref{asymp-unc-part4}
with a stronger error term of the form 
$N \exp ( - C  (\log N) ^{3/5}(\log \log N)^{-1/5})$.
As in the previous case,
the comparative weakness of our error term is due to the
use of the zero-density estimates for the Riemann zeta-function.
A direct trial following the lines of  Lemma 11 of  Plaksin \cite{Plaksin1981}
leads to weaker uniformity ranges:  
$H\gg N^{3/4}L^A$, for some $A>0$, assuming RH and  $H\gg N^{7/24+1/2+\eps}$ unconditionally. 

Concerning the  problem about the sum of two squares,
i.e. the asymptotic formula for
\[
r_{2,2}(n) =
\sum_{m_{1}^{2}+m_{2}^2=n} 
1,
\]
our method leads to a weaker result than the one that follows from 
the well-known formula 
\(
\sum_{n=1}^{N}  \, r_{2,2}(n)
=
\frac{\pi}{4} N - N^{1/2} + \Odi{ N^{\alpha} },
\)
with $\alpha\in (1/4,1/3)$.

In the proofs we will use the original Hardy-Littlewood circle method setting. 
This depends on the fact in the standard finite sums method the approximation 
needed to detect the main term contribution leads to an error term which is under
control essentially only for $H>N^{2/3+\eps}$, see also Remark \ref{Remark-dopo-teo2} at the bottom of the proof of Theorem \ref{asymp-unc-part3}.   

\medskip
{\bf Acknowledgements.}    This research was partially supported by the grant PRIN2010-11 \textsl{Arithmetic Algebraic Geometry and Number Theory}. We wish to thank the referee for his/her remarks.

\section{Definitions and Lemmas}  

Let $L=\log N$ and $\ell\ge 1$ be an integer. We define  
\begin{align} 
\label{main-defs}
\Stilde_\ell(\alpha) = \sum_{n=1}^{\infty} \Lambda(n) e^{-n^{\ell}/N} e(n^{\ell}\alpha) , & \quad
z= 1/N-2\pi i\alpha
\quad
\textrm{and}
\quad 
U(\alpha,H) = \sum_{1\leq m\leq H}e(m\alpha),
\end{align}
where  $e(\alpha) = e^{2\pi i\alpha}$. {}From now on, we denote 
\begin{equation}
\label{E-def}
\Etilde_{\ell}(\alpha) : =\Stilde_\ell(\alpha) - \frac{\Gamma(1/\ell)}{\ell z^{1/\ell}}.
\end{equation}

We will also need the following unconditional version of Lemma 3
of \cite{LanguascoZ2016a}; the proof is essentially the same used there and so
we skip part of the argument.  We just repeat the definition
of the main quantities involved and  write how to use the zero-density estimates 
to conclude the proof.
In fact all of the following lemmas will be used  just for $\ell=1,2$ but we take this occasion
to describe the  general case.
\begin{Lemma}
\label{LP-Lemma-gen-unc} 
Let $\eps$ be an arbitrarily small
positive constant,  $\ell \ge 1$ be an integer and $N$ be a
sufficiently large integer. Then there exists a positive constant 
$c_1 = c_{1}(\eps)$, which does not depend on $\ell$, such that 
\[
\int_{-\xi}^{\xi} \,
 \vert
\Etilde_{\ell}(\alpha) 
 \vert^{2}\,
\dx \alpha 
\ll_{\ell}
 N^{2/\ell-1} \exp \Big( - c_{1}  \Big( \frac{L}{\log L} \Big)^{1/3} \Big) 
\]
uniformly for $ 0\le \xi < N^{-1 +5/(6\ell) - \eps}$.
\end{Lemma} 

\begin{Proof}
Since 
$z^{-\rho/\ell} =  \vert z \vert^{-\rho/\ell} \exp\bigl(-i(\rho/\ell)\arctan2\pi N\alpha\bigr)$,
by  Stirling's formula we have that
\[
\frac{1}{\ell}
\sum_{\rho}z^{-\rho/\ell}\Gamma\Bigl(\frac{\rho}{\ell}\Bigr) 
\ll_{\ell}
\sum_{\rho}  \vert z \vert^{-\beta/\ell}
\vert \gamma \vert ^{\beta/\ell -1/2}
\exp
\Bigl(\frac{\gamma}{\ell}\arctan2\pi N\alpha - \frac{\pi}{2\ell} \vert\gamma \vert\Bigr).
\]
Recalling  the Vinogradov-Korobov zero-free region, \emph{i.e.},
there are no zeros $\beta+i\gamma$ of the Riemann zeta function having 
\begin{equation}
\label{Vino-Koro}
\beta 
\geq 
1- \frac{c^{\prime}}{(\log ( \vert \gamma \vert +2))^{2/3}(\log \log ( \vert \gamma \vert +2))^{1/3}}
= 1-\delta(\gamma),
\end{equation} 
say, where $c^{\prime}>0$ is an absolute constant, for
$\vert\alpha \vert\leq\ 1/N$ or $\gamma \alpha <0$ we get
\begin{align*}
\sum_{\rho}z^{-\rho/\ell}\Gamma (\rho/\ell) 
&\ll
N^{1/\ell} 
\sum_{\rho}  N^{-\delta(\gamma)/\ell}
\vert \gamma \vert ^{1/\ell - 1/2} 
\exp
\Bigl(-C\frac{\vert\gamma \vert}{\ell}\Bigr)
\\
& \ll_{\ell}
N^{1/\ell}  
\sum_{\rho}  N^{-\delta(\gamma)/\ell}
\exp \Bigl(-C_{1}\frac{\vert\gamma \vert}{\ell}\Bigr)
\ll_{\ell} 
N^{(1-\eps)/\ell} ,
\end{align*}
where
$C,C_{1}>0$ are  absolute positive constants and $\eps \in (0,1)$ is suitably small. 
Hence, by the explicit formula for $\Stilde_\ell$ which is Lemma 2 of
\cite{LanguascoZ2016a}, we have
\begin{equation}
\label{R1-first-unc} 
I(N,\xi,\ell)
:=
\int_{-\xi}^{\xi} \,
\Bigl\vert
\Stilde_\ell(\alpha) - \frac{\Gamma(1/\ell)}{\ell z^{1/\ell}}
\Bigr\vert^{2} \dx \alpha 
\ll_{\ell} 
N^{2(1-\eps)/\ell} \xi
\end{equation}
if $0\leq\xi\le 1/N$, and
\begin{equation}
\label{R1-estim-unc}
I(N,\xi,\ell)
\ll_{\ell}
\int_{1/N}^{\xi}
\Big \vert\sum_{\rho\colon \gamma > 0}z^{-\rho/\ell}\Gamma\Bigl(\frac{\rho}{\ell}\Bigr) \Big \vert^2 \dx \alpha +
\int_{-\xi}^{-1/N}
\Big \vert\sum_{\rho\colon \gamma < 0}z^{-\rho/\ell}\Gamma\Bigl(\frac{\rho}{\ell}\Bigr) \Big \vert^2 \dx \alpha 
+ 
N^{2/\ell - 1 -2 \eps/\ell}
 \end{equation}
if $\xi>1/N$. We will treat only the first integral on the right hand
side of \eqref{R1-estim-unc}, the second being completely similar. 
Clearly
\begin{equation}
\label{split-int-unc} 
\int_{1/N}^{\xi}
\Big \vert\sum_{\rho\colon \gamma > 0}z^{-\rho/\ell}\Gamma\Bigl(\frac{\rho}{\ell}\Bigr) \Big \vert^2 
\dx \alpha 
= 
\sum_{k=1}^K
\int_\eta^{2\eta} \Big \vert\sum_{\rho\colon \gamma > 0}z^{-\rho/\ell}\Gamma\Bigl(\frac{\rho}{\ell}\Bigr) \Big \vert^2 \dx \alpha
 +
\Odi{1} 
 \end{equation}
where $\eta=\eta_k= \xi/2^k$, $1/N\le \eta \le \xi/2$  and $K$ is a suitable integer satisfying $K=\Odi{L}$. 
Writing $\arctan 2\pi N\alpha = \pi/2 - \arctan(1/2\pi N\alpha)$ and
using the Saffari-Vaughan technique we have
\begin{align}
\notag
\int_{\eta}^{2\eta} 
\Big \vert\sum_{\rho\colon \gamma > 0}z^{-\rho/\ell}\Gamma\Bigl(\frac{\rho}{\ell}\Bigr) 
\Big \vert^2\ \dx \alpha 
&\le 
\int_1^2 
\Bigl(
\int_{\delta\eta/2}^{2\delta\eta}
\Big \vert\sum_{\rho\colon \gamma > 0}z^{-\rho/\ell}\Gamma\Bigl(\frac{\rho}{\ell}\Bigr) \Big \vert^2 \dx \alpha 
\Bigr) 
\dx \delta  
 \\
\label{expl-unc} 
&= \sum_{\rho_1\colon \gamma_{1} > 0}\sum_{\rho_2\colon \gamma_{2} > 0}
\Gamma\Bigl(\frac{\rho_{1}}{\ell}\Bigr)
\overline{\Gamma\Bigl(\frac{\rho_{2}}{\ell}\Bigr)}\ 
e^{\frac{\pi}{2\ell}(\gamma_1+\gamma_2 -i (\beta_{1}-\beta_{2}))}\, \cdot  J,
\end{align}
say, where 
\[
J = J(N,\eta,\ell,\beta_1,\beta_2,\gamma_1,\gamma_2) 
= 
\int_1^2
\Bigl(
\int_{\delta\eta/2}^{2\delta\eta}
f_{1}(\alpha)f_2(\alpha)\ \dx \alpha 
\Bigr)
\dx \delta,
\]
\[
f_{1}(\alpha)=  
\vert z \vert^{-w} ,
\quad 
f_{2}(\alpha) = 
\exp\Bigl(-\frac{\gamma_1+\gamma_2 -i (\beta_{1}-\beta_{2})}{\ell}\arctan\frac{1}{2\pi N\alpha}\Bigr),
\]
and
$w=w(\ell,\beta_1,\beta_2,\gamma_1,\gamma_2)=
(\beta_{1}+\beta_{2})/\ell+(i/\ell)(\gamma_1-\gamma_2)$.
Arguing exactly as in the proof of Lemma 3 in \cite{LanguascoZ2016a}, see pages 6-7 there, 
we get
\[
J 
\ll_{\ell}
\eta^{1-(\beta_{1}+\beta_{2})/\ell} 
\frac{1+\ (\frac{1+ \gamma_1+\gamma_2}{N\eta})^2} {1+ \vert\gamma_1-\gamma_2 \vert^2}
\exp\Bigl(-c\Bigl(\frac{ \gamma_1  +  \gamma_2}{N\eta}\Bigr)\Bigr),
\]
hence from \eqref{expl-unc}  and Stirling's formula we have
\begin{align}
\notag
\int_{\eta}^{2\eta} 
\Big \vert
\sum_{\rho \colon \gamma>0}z^{-\rho/\ell}\Gamma\Bigl(\frac{\rho}{\ell}\Bigr) 
\Big \vert^2
\dx \alpha 
\ll_{\ell}   
\sum_{\rho_{1}\colon \gamma_{1} > 0} &\sum_{\rho_{2}\colon \gamma_{2} > 0}
\eta^{1-(\beta_{1}+\beta_{2})/\ell} 
 \gamma_1 ^{\beta_{1}/\ell-1/2}  \gamma_2 ^{\beta_{2}/\ell-1/2}
\\
\label{int-estim-unc} 
&\times
\frac{1+(\frac{1+ \gamma_1+\gamma_2}{N\eta})^2}{1+ \vert\gamma_1-\gamma_2 \vert^2}
\exp\Bigl(-c\Bigl(\frac{ \gamma_1  +  \gamma_2}{N\eta}\Bigr)\Bigr) .
\end{align}
Sorting real and imaginary parts it is clear that
\[
 \gamma_1 ^{\beta_{1}/\ell-1/2} 
 \gamma_2 ^{\beta_{2}/\ell-1/2}
\Bigl\{
1+\Bigl(\frac{1+ \gamma_1+\gamma_2}{N\eta}\Bigr)^2
\Bigr\}
\exp\Bigl(-c\Bigl(\frac{ \gamma_1  +  \gamma_2}{N\eta}\Bigr)\Bigr) 
\ll_{\ell}
  \gamma_1 ^{2\beta_{1}/\ell -1}
\exp\Bigl(-\frac{c}{2}\frac{ \gamma_1 }{N\eta}\Bigr),
\]
hence the r.h.s.~of \eqref{int-estim-unc}  becomes
\begin{align}
\notag
\ll_{\ell} 
&
\sum_{\rho_{1}\colon \gamma_{1}>0}
\eta^{1-2\beta_{1}/\ell} 
  \gamma_1 ^{2\beta_{1}/\ell -1}
 \exp\Bigl(-\frac{c}{2}\frac{ \gamma_1 }{N\eta}\Bigr)
\sum_{\rho_{2}\colon \gamma_{2}>0; \beta_{2}\le \beta_{1}}
\frac{1}{1+ \vert\gamma_1-\gamma_2 \vert^2} 
\\
\label{int-estim1-unc} 
&
\ll_{\ell} 
\sum_{\rho_{1}\colon \gamma_{1}>0}
\Bigr(\frac{ \gamma_1 }{\eta}\Bigr)^{2\beta_{1}/\ell -1}
 \exp\Bigl(-\frac{c}{4}\frac{ \gamma_1 }{N\eta}\Bigr)
\end{align}
since the number of zeros $\rho_2= \beta_{2}+i\gamma_2$ with $n\le  \vert\gamma_1-\gamma_2 \vert\le n+1$
is $\Odi{\log (n+ \gamma_1 )}$.

Now we  use \eqref{Vino-Koro} and  the Ingham-Huxley zero-density estimate, \emph{i.e.}, for  $1/2\leq \sigma \leq 1$ we have
that $N(\sigma,t) \ll t^{(12/5)(1-\sigma)} (\log t)^{B}$.
Hence, uniformly for $1/N < \eta < N^{-1 +5/(6\ell) - \eps}$, by 
(6.17) of Saffari and Vaughan \cite{SaffariV1977}
we get that \eqref{int-estim1-unc} is
\begin{align}
 \notag
\sum_{\rho_{1}\colon \gamma_{1}>0} 
&
\Bigr(\frac{ \gamma_1 }{\eta}\Bigr)^{2\beta_{1}/\ell -1} 
 \exp\Bigl(-\frac{c}{4}\frac{ \gamma_1 }{N\eta}\Bigr)
  \ll
\sum_{\substack{\beta_{1}\ge 1/2\\ 0< \gamma_1  \le N^{4}} }
\Bigr(\frac{ \gamma_1 }{\eta}\Bigr)^{2\beta_{1}/\ell -1} 
 \exp\Bigl(-\frac{c}{4}\frac{ \gamma_1 }{N\eta}\Bigr)
\\
 \notag
& \ll_{\ell}
\max_{1/2\le \sigma \le 1-\delta(N^{4})} 
\int_{0}^{N^{4}} 
t^{(12/5)(1-\sigma)} (\log t)^{B}
\Bigl[
\Bigr(\frac{t}{\eta}\Bigr)^{2\sigma/\ell -1} 
 \exp\Bigl(-\frac{c}{8}\frac{t}{N\eta}\Bigr)
 \Bigr]^{\prime}
\ \dx t
 \\
  \notag
 & \ll_{\ell} 
\max_{1/2\le \sigma \le 1-\delta(N^{4})} 
\int_{0}^{\infty} 
(Nu)^{2\sigma/\ell -1} (N\eta u)^{(12/5)(1-\sigma)}
 \exp\Bigl(-\frac{c}{8}u\Bigr)
\ \dx u
 \\
\label{int-estim2-unc}  
& \ll_{\ell} 
\max_{1/2\le \sigma \le 1-\delta(N^{4})}  
\Bigl((N\eta)^{(12/5)(1-\sigma)} N^{2\sigma/\ell -1} \Bigr)
 \ll_{\ell}
 N^{2/\ell-1} \exp \Bigl( - c_{1}  \Big( \frac{L}{\log L} \Big)^{1/3} \Bigr),
\end{align}
where $c_1 = c_{1}(\eps)$ is a positive constant which does not depend on $\ell$.
{}From \eqref{R1-first-unc}-\eqref{split-int-unc}  
and \eqref{int-estim-unc}-\eqref{int-estim2-unc}  we get
\begin{equation}
\label{last-unc} 
\int_{-\xi}^{\xi} 
\Big \vert
\sum_{\rho\colon \gamma>0}z^{-\rho/\ell}\Gamma\Bigl(\frac{\rho}{\ell}\Bigr) 
\Big \vert^2  
\dx \alpha 
\ll_{\ell} 
 N^{2/\ell-1} \exp \Big( - c_{1}  \Big( \frac{L}{\log L} \Big)^{1/3} \Big)
\end{equation}
uniformly for  $1/N< \xi < N^{-1 +5/(6\ell) - \eps}$. 
Lemma \ref{LP-Lemma-gen-unc}   follows from \eqref{R1-first-unc}-\eqref{R1-estim-unc} and \eqref{last-unc}.
\end{Proof}

We need also  the following  analogue of Lemma 1 
of \cite{LanguascoZ2015c}. 
Let
\begin{equation}
\label{omega-def}
\omega_{\ell}(\alpha)
=
\sum_{m=1}^{\infty}   
e^{-m^{\ell}/N} e(m^{\ell}\alpha)
=
\sum_{m=1}^{\infty}   
e^{-m^{\ell}z} .
\end{equation}
We explicitly remark that for $\ell=1$ the  proof of 
Lemma \ref{zac-lemma-series} gives just trivial results; in this case
a non-trivial estimate, which, in any case, is not useful in this context, 
can be obtained following the line of Corollary 3 of \cite{LanguascoP1994}.
\begin{Lemma}
\label{zac-lemma-series}
Let $\ell\ge 2 $ be an integer and $0<\xi\leq 1/2$. Then
\[
\int_{-\xi}^{\xi} 
|\omega_{\ell}(\alpha)|^2 \ \dx\alpha 
\ll_{\ell}
  \xi N^{1/\ell}
  + 
\begin{Biggcases}
L & \text{if}\ \ell =2\\
1 & \text{if}\ \ell > 2
\end{Biggcases}\]
and
\[
\int_{-\xi}^{\xi} 
|\Stilde_{\ell}(\alpha)|^2 \ \dx\alpha 
\ll_{\ell}
\xi N^{1/\ell} L  +
\begin{Biggcases}
L^{2} & \text{if}\ \ell =2\\
1 & \text{if}\ \ell > 2.
\end{Biggcases}
\]
\end{Lemma}
\begin{Proof}
By symmetry we can  integrate over $[0,\xi]$.
We use Corollary 2 of Montgomery and Vaughan  \cite{MontgomeryV1974} 
(see also the remark after their statement) 
with $T=\xi$, $a_r=\exp(-r^\ell/N)$ and $\lambda_r= 2\pi r^\ell$ thus getting
\begin{align*}
\int_{0}^{\xi} \vert \omega_{\ell}(\alpha) \vert ^2\, \dx \alpha 
&=
\sum_{r\ge 1} e^{-2r^\ell/N} \bigl(\xi +\Odim{\delta_r^{-1}}\bigr)
\ll_\ell
\xi N^{1/\ell} 
+
\sum_{r\ge 1} r^{1-\ell} e^{-2r^\ell/N}
\end{align*}
since $\delta_{r} = \lambda_r - \lambda_{r-1}  \gg_{\ell} r^{\ell-1}$. 
The last  term is $\ll_{\ell}1$ if $\ell >2$ and $\ll L$ otherwise.
This proves the first part of Lemma \ref{zac-lemma-series}.
Arguing analogously with  $a_r=\Lambda(r)\exp(-r^\ell/N)$, 
by the Prime Number Theorem we get
\begin{align*}
\int_{0}^{\xi} \vert \Stilde_{\ell}(\alpha) \vert ^2\, \dx \alpha 
&=
\sum_{r\ge 1} \Lambda(r)^2 e^{-2r^\ell/N} \bigl(\xi +\Odim{\delta_r^{-1}}\bigr)
\ll_\ell
\xi N^{1/\ell} L
+
\sum_{r\ge 1} \Lambda(r)^2 r^{1-\ell} e^{-2r^\ell/N}.
\end{align*}
The last  term is $\ll_{\ell}1$ if $\ell >2$ and $\ll L^{2}$ otherwise.
The second part of Lemma \ref{zac-lemma-series} follows.
\end{Proof}
Let now
\begin{equation}
\label{T-tilde-def}
\Ttilde_\ell(\alpha) = \sum_{p=2}^{\infty} \log p \, e^{-p^{\ell}/N} e(p^{\ell}\alpha).
\end{equation}
We also have
 \begin{Lemma}
\label{trivial-lemma}
Let $\ell\ge 1 $ be an integer. Then
\(
\vert \Stilde_{\ell}(\alpha)- \Ttilde_{\ell}(\alpha) \vert 
\ll_{\ell}
 N^{1/(2\ell)}  .
\)
\end{Lemma}
\begin{Proof}
Clearly we have
\begin{align*}
 \vert \Stilde_{\ell}(\alpha)- \Ttilde_{\ell}(\alpha) \vert 
\le 
\sum_{k\ge 2}\sum_{p\ge 2}   \log p\,  e^{-p^{k\ell}/N}   
\ll_\ell
 N^{1/(2\ell)}  
 \end{align*}
where in the last inequality we used  the Prime Number Theorem.
\end{Proof}

Letting $\omega(\alpha)=\omega_{2}(\alpha)$ and 
\[
\theta(z)
=
\sum_{n=-\infty}^{\infty}   
e^{-n^{2}/N} e(n^{2}\alpha)
=
\sum_{n=-\infty}^{\infty}   
e^{-n^{2}z} 
=
1+2\omega(\alpha),
\]
the functional equation of the $\theta$-function 
(see, \emph{e.g.}, Proposition VI.4.3, page 340, of Freitag and Busam
\cite{FreitagB2009}) gives  that
\(
\theta(z)
=
(\pi/z)^{1/2} 
\theta ( \pi^2/z).
\)
Hence  we have
\begin{equation}
\label{omega-approx}
\omega(\alpha)  
= 
\frac12
\left(\frac{\pi}{z}\right)^{1/2} 
\!\!\!\!
- \frac12
+
\left(\frac{\pi}{z}\right)^{1/2} 
\sum_{\ell=1}^{+\infty}   
e^{-\ell^{2}\pi^{2}/z}.
\end{equation} 
\begin{Lemma} 
\label{omega-Y}
Let $N$ be a large integer,  $z= 1/N-2\pi i\alpha$, $\alpha\in [1/2,1/2]$ and 
$Y=\Re(1/z)>0$.
We have 
\[
\Bigl\vert 
\sum_{\ell=1}^{+\infty}   
e^{-\ell^{2}\pi^{2}/z}
\Bigr\vert  
\ll 
\begin{Biggcases} 
e^{- \pi^{2}  Y } & \textrm{for} \ Y\ge 1 \\
Y^{-1/2} & \textrm{for} \ 0 <Y\le 1.
\end{Biggcases}
\]
\end{Lemma} 
\begin{Proof}
It is clear that
\begin{align*}
\Bigl\vert 
\sum_{\ell=1}^{+\infty}   
e^{-\ell^{2}\pi^{2}/z}
\Bigr\vert  
&\le 
\sum_{\ell=1}^{+\infty}   
e^{-\ell^{2}\pi^{2} Y}
\le
\sum_{\ell=1}^{+\infty}   
e^{-\ell\pi^{2} Y}
=
\frac{e^{- \pi^{2}  Y } }{1-e^{- \pi^{2}  Y } }
\ll
e^{-\pi^{2}  Y }
\end{align*}
for $Y\ge 1$. Moreover, for $Y>0$, we also have
\[
\sum_{\ell=1}^{+\infty}   
e^{-\ell^{2}\pi^{2} Y}
\le
1+ \int_{1}^{+\infty}  e^{-t^{2}\pi^{2} Y} \dx t
\ll
1+  Y^{-1/2}
\]
and the lemma is proved.
\end{Proof}

Since
\begin{equation}
\notag
Y
=
\Re(1/z)
=
\frac{N}{1+4\pi^2\alpha^2N^2}
\ge
\frac{1}{5\pi^2}
\begin{Biggcases}
N & \textrm{if} \ \vert \alpha \vert \le 1/N\\
(\alpha^2N)^{-1} & \textrm{if} \ \vert \alpha \vert > 1/N,
\end{Biggcases}
\end{equation}
from Lemma \ref{omega-Y} we get 
\begin{equation}
\label{omega-Y-estim}
\Bigl\vert 
\sum_{\ell=1}^{+\infty}   
e^{-\ell^{2}\pi^{2}/z}
\Bigr\vert  
\ll
\begin{cases}
\exp(-\pi^2 N) & \textrm{if} \ \vert \alpha \vert \le 1/N\\
\exp(-\pi^2/(\alpha^2N)) & \textrm{if} \ 1/N<\vert \alpha \vert = \odim{N^{-1/2}}\\
1+ N^{1/2} \vert \alpha\vert & \textrm{otherwise}.
\end{cases}
\end{equation}
We also recall that
\begin{equation}
\label{UH-estim}
\vert U(\alpha,H) \vert
\le
\min \bigl(H; \vert \alpha \vert^{-1} \bigr),
\end{equation}
\begin{equation}
\label{z-estim}
\vert z\vert ^{-1} \ll \min \bigl(N, \vert \alpha \vert^{-1}\bigr)
\end{equation}
and we finally define
\begin{equation}
\label{B-def}
B=B(N,c)= \exp \Big( c   \Big( \frac{L}{\log L} \Big)^{1/3} \Big),
\end{equation}
where $c=c(\eps)>0$ will be chosen later.

\section{Proof of Theorem \ref{asymp-part3} }  

%
Recalling \eqref{main-defs} and \eqref{T-tilde-def},
it is an easy matter to see that
\begin{align}
\notag
\sum_{n = N+1}^{N + H} & e^{-n/N}
r^{\second}_{2,2}(n) 
= 
\int_{-1/2}^{1/2}  \Ttilde_2(\alpha)^{2}  U(-\alpha,H) e(-N\alpha) \, \dx \alpha
\\
\notag
 &=  \int_{-1/2}^{1/2} ( \Ttilde_2(\alpha)^{2}- \Stilde_2(\alpha)^{2})  U(-\alpha,H) e(-N\alpha)\, \dx \alpha
\\
\notag
 & \hskip1cm +
 \int_{-1/2}^{1/2}      \frac{\pi}{4z} U(-\alpha,H)
 e(-N\alpha)\, \dx \alpha
  +
 \int_{-1/2}^{1/2} \Bigl( \Stilde_{2}(\alpha)^{2} -  \frac{\pi}{4z}	\Bigr) U(-\alpha,H)
 e(-N\alpha)\, \dx \alpha 
\\
\label{approx-th1-part3} 
&
= I_{0} + I_{1}+I_{2}  , 
\end{align}
say.
Using  the identity $f^{2}-g^{2} = 2f(f-g) -(f-g)^{2}$ and the Cauchy-Schwarz inequality we have
\begin{align*} 
I_{0}
&\ll
 \int_{-1/2}^{1/2}
      \vert \Stilde_{2}(\alpha)  \vert \vert \Stilde_{2}(\alpha) -\Ttilde_{2}(\alpha)  \vert
        \vert U(\alpha,H)\vert 
     \, \dx \alpha
+
 \int_{-1/2}^{1/2}
       \vert \Stilde_{2}(\alpha) -\Ttilde_{2}(\alpha)  \vert^{2}
       \vert U(\alpha,H) \vert
        \, \dx \alpha 
\\  
& \ll
N^{1/4}
\Bigl(
\int_{-1/2}^{1/2}
      \vert \Stilde_{2}(\alpha)  \vert^{2}\vert U(\alpha,H) \vert
        \dx \alpha
     \Bigr)^{1/2}
 \Bigl(
      \int_{-1/2}^{1/2}
       \vert U(\alpha,H) \vert
        \, \dx \alpha        
\Bigr)^{1/2}
\!\!\! 
+
N^{1/2}
 \int_{-1/2}^{1/2}
       \vert U(\alpha,H) \vert
        \, \dx \alpha,
\end{align*}
by Lemma \ref{trivial-lemma}.
By Lemma \ref{zac-lemma-series}, \eqref{UH-estim} and
a partial integration argument we obtain
\begin{align*}
\int_{-1/2}^{1/2}
\vert \Stilde_{2}(\alpha)  \vert^{2}\vert U(\alpha,H) \vert
        \dx \alpha
&\ll
H \int_{-1/H}^{1/H}
\vert \Stilde_{2}(\alpha)  \vert^{2}\vert U(\alpha,H) \vert
        \dx \alpha
 +
 \int_{1/H}^{1/2}
\vert \Stilde_{2}(\alpha)  \vert^{2}\frac{\dx \alpha}{\alpha}
\\
&\ll
H\Bigl(\frac{N^{1/2}L}{H}+L^{2}\Bigr) 
+ 
 N^{1/2}L 
+
 \int_{1/H}^{1/2}
 (N^{1/2} \xi L+L^{2}) \frac{\dx \xi}{\xi^{2}}
 \\
&\ll  
 N^{1/2}L^{2} + HL^{2}.
\end{align*}
Hence
\begin{equation}
\label{I0-estim}
I_{0}
\ll
N^{1/4}(N^{1/2}L^2 +H L^{2})^{1/2}L^{1/2} + N^{1/2}L
\ll
N^{1/2}L^{3/2}+ H^{1/2}N^{1/4}L^{3/2}.
\end{equation}

Now we evaluate $I_{1}$. Using Lemma~4 of \cite{LanguascoZ2016a} 
we immediately get 
\begin{equation}
\label{I1-eval-part3}
I_{1}
=
\frac{\pi}{4}   \sum_{n = N+1}^{N + H}  e^{-n/N} + \Odi{\frac{H}{N}}
  =
\frac{\pi H }{4e} + \Odi{\frac{H^2}{N}}.
\end{equation} 

Now we estimate $I_{2}$.
Again using  the identity $f^{2}-g^{2} = 2f(f-g) -(f-g)^{2}$, by \eqref{E-def} we obtain 
\begin{equation}
\label{I2-estim1-part3}
I_{2}
 \ll
 \int_{-1/2}^{1/2}
      \vert \Etilde_{2}(\alpha)  \vert
        \frac{\vert U(\alpha,H)\vert}{\vert z\vert^{1/2}}
     \, \dx \alpha
+
 \int_{-1/2}^{1/2}
       \vert \Etilde_{2}(\alpha)  \vert^{2}  
       \vert U(\alpha,H) \vert
        \, \dx \alpha 
=
   J_{1}+J_{2},
\end{equation}
say.
Using \eqref{UH-estim}-\eqref{z-estim}, Lemma 3 of \cite{LanguascoZ2016a} and a partial integration argument we have
\begin{equation}
 \label{J2-estim-part3}
J_{2}
\ll 
H
 \int_{-1/H}^{1/H}  \vert \Etilde_{2}(\alpha)  \vert^{2}   \, \dx \alpha 
 +
 \int_{1/H}^{1/2}  \vert \Etilde_{2}(\alpha)  \vert^{2}   
 \frac{\dx \alpha}{\alpha}
\ll 
N^{1/2} L^{2} 
+
N^{1/2} L^{2} \Bigl(1 +   \int_{1/H}^{1/2}   \frac{\dx \xi}{\xi} \Bigr)
 \ll 
 N^{1/2} L^{3}. 
\end{equation}
Using  the Cauchy-Schwarz inequality and arguing as for $J_{2}$ we get 
\begin{align}
\notag
J_{1}
&\ll 
H N^{1/2} \Bigl(  \int_{-1/N}^{1/N}\!\!\!\! \dx \alpha  \Bigr)^{1/2}
\Bigl(
 \int_{-1/N}^{1/N} \!\! \vert \Etilde_{2}(\alpha)  \vert^{2}   \, \dx \alpha 
 \Bigr)^{1/2}
 \!\!\!\!\! +
 H \Bigl(  \int_{1/N}^{1/H} \!\!  \frac{\dx \alpha}{\alpha^{1/2}}\Bigr)^{1/2} 
\Bigl(
 \int_{1/N}^{1/H}\!\!  \vert \Etilde_{2}(\alpha)  \vert^{2}  
 \frac{\dx \alpha}{\alpha^{1/2}}
  \Bigr)^{1/2}
  \\
  \notag
     &\hskip2cm 
     +
  \Bigl(  \int_{1/H}^{1/2}   \frac{\dx \alpha}{\alpha^{3/2}}\Bigr)^{1/2} 
\Bigl(
 \int_{1/H}^{1/2}  \vert \Etilde_{2}(\alpha)  \vert^{2}   
 \frac{\dx \alpha}{\alpha^{3/2}}
   \Bigr)^{1/2}
 \\
 \notag 
 &\ll 
H N^{-1/4} L 
+ 
H^{3/4}  N^{1/4}  L  \Bigl(  \frac{1}{H}+ \int_{1/N}^{1/H}  \frac{\dx \xi}{\xi^{1/2}}   \Bigr)^{1/2} 
 +
H^{1/4} N^{1/4} L \Bigl( H^{1/2} +   \int_{1/H}^{1/2}  \frac{\dx \xi}{\xi^{3/2}}   \Bigr)^{1/2}
 \\
 \label{J1-estim-part3}
 &\ll  
 H^{1/2}N^{1/4} L.
\end{align}

Combining \eqref{I2-estim1-part3}-\eqref{J1-estim-part3} we finally obtain
\begin{equation}
\label{I2-estim-final-part3}
I_{2} \ll  H^{1/2}N^{1/4} L + N^{1/2} L^{3} .
\end{equation}
Now using \eqref{approx-th1-part3}-\eqref{I1-eval-part3}
and \eqref{I2-estim-final-part3} we  have 
\begin{equation} 
\label{almost-done}
\sum_{n = N+1}^{N + H}  e^{-n/N}
r^{\second}_{2,2}(n)   
 =
\frac{\pi H }{4e} + \Odig{\frac{H^2}{N}  + N^{1/2} L^{3} +  H^{1/2}N^{1/4} L^{3/2}}
\end{equation}
which is an asymptotic formula for $\infty(N^{1/2}L^{3}) \le H \le \odi{N}$.
{}From  $e^{-n/N}=e^{-1}+ \Odi{H/N}$ for $n\in[N+1,N+H]$,
we get 
\begin{equation}
\label{th1-final} 
   \sum_{n = N+1}^{N + H} 
r^{\second}_{2,2}(n) 
    = 
    \frac{\pi H }{4} 
    +
    \Odig{\frac{H^2}{N}  + N^{1/2} L^{3} +  H^{1/2}N^{1/4} L^{3/2}}
  +
  \Odig{\frac{H}{N}\sum_{n = N+1}^{N + H} r^{\second}_{2,2}(n)}.
\end{equation}
Using $e^{n/N}\leq e^{2}$ 
and \eqref{almost-done} for $H$ in the previously mentioned range,
it is easy to see that the last error term is
$\ll H^{2}N^{-1}$.
Combining \eqref{th1-final} and the last remark,
Theorem \ref{asymp-part3}  hence follows for $\infty(N^{1/2}L^{3}) \le H \le \odi{N}$.
%
\qed

\section{Proof of Theorem \ref{asymp-unc-part3}}  
 
Recalling \eqref{main-defs} and \eqref{T-tilde-def}, it is an easy matter to see that
\begin{align}
\notag
\sum_{n = N+1}^{N + H}&  e^{-n/N}
r^{\second}_{2,2}(n) 
= 
\int_{-1/2}^{1/2}  \Ttilde_2(\alpha)^{2}  U(-\alpha,H) e(-N\alpha) \, \dx \alpha
\\
\notag
 &=
\int_{-1/2}^{1/2} (\Ttilde_2(\alpha)^{2}-  \Stilde_2(\alpha)^{2})  U(-\alpha,H) e(-N\alpha)\, \dx \alpha  
+
 \int_{-B/H}^{B/H}     \frac{\pi}{4z}  U(-\alpha,H)
  e(-N\alpha)\, \dx \alpha
\\
 \notag
 & \hskip1cm
  +
 \int_{-B/H}^{B/H}
 \Bigl( \Stilde_{2}(\alpha)^{2} -  \frac{\pi}{4z}	\Bigr) U(-\alpha,H)
 e(-N\alpha)\, \dx \alpha 
+ 
\int\limits_{\mathclap{[-1/2,-B/H]\cup  [B/H,1/2]}} 
\Stilde_2(\alpha)^{2} U(-\alpha,H) e(-N\alpha) \, \dx \alpha
\\
\label{approx-th2-part3} 
&
= I_{0}+I_{1}+I_{2}+I_{3}  , 
\end{align}
say, where $B$ is defined in \eqref{B-def}.
$I_{0}$ can be estimated as in \eqref{I0-estim} and gives
\begin{equation}
\label{I0-estim-unc}
I_{0}
\ll
N^{1/2}L^{3/2} + H^{1/2}N^{1/4}L^{3/2}.
\end{equation}

Now we evaluate $I_{1}$. Using Lemma~4 of \cite{LanguascoZ2016a} 
and \eqref{UH-estim}
we immediately get 
\begin{equation}
\label{I1-eval-part3-unc}
I_{1}
=
\frac{\pi}{4}   \sum_{n = N+1}^{N + H}    e^{-n/N} + \Odig{\frac{H}{N}}
+
\Odig{\int_{B/H}^{1/2} \frac{\dx\alpha}{\alpha^2}}
  =
\frac{\pi H}{4e} + \Odig{\frac{H^2}{N}+\frac{H}{B}}.
\end{equation}

Now we estimate $I_{2}$.
Using the identity $f^{2}-g^{2} = 2f(f-g) -(f-g)^{2}$, by \eqref{E-def} and  \eqref{UH-estim} we obtain 
\begin{equation}
\label{I2-estim1-part3-unc}
I_{2}
\ll
H
\Bigl(
 \int_{-B/H}^{B/H}
      \vert \Etilde_{2}(\alpha)  \vert
        \frac{\dx \alpha}{\vert z\vert^{1/2}}
+ 
 \int_{-B/H}^{B/H}
       \vert \Etilde_{2}(\alpha)  \vert^{2}  
        \, \dx \alpha 
\Bigr)
   =
  H( J_{1}+J_{2}),
\end{equation}
say.
Using  Lemma \ref{LP-Lemma-gen-unc} with $\ell=2$  we have
\begin{align}
 \label{J2-estim-part3-unc}
J_{2} 
\ll
 \exp \Big( - c_{1}  \Big( \frac{L}{\log L} \Big)^{1/3} \Big) 
\end{align}
provided that $\infty(1/N)<B/H < N^{-7/12 - \eps/2}$, i.e. 
$N^{7/12 + \eps} \le H \le \odi{N}$ suffices.

Using  the Cauchy-Schwarz inequality and arguing as for $J_{2}$ we get 
\begin{equation}
 \label{J1-estim-part3-unc}
J_{1}
\ll  \Bigl(  \int_{-B/H}^{B/H} \frac{\dx \alpha}{\vert z \vert}  \Bigr)^{1/2}
\Bigl(
 \int_{-B/H}^{B/H} \!\! \vert \Etilde_{2}(\alpha)  \vert^{2}   \, \dx \alpha 
 \Bigr)^{1/2} 
 \ll  
 \exp \Big( - \frac{c_{1}}{4}  \Big( \frac{L}{\log L} \Big)^{1/3} \Big),
\end{equation}
provided that $\infty(1/N)<B/H < N^{-7/12 - \eps/2}$, i.e. 
$N^{7/12 + \eps} \le H \le \odi{N}$.

Combining \eqref{I2-estim1-part3-unc}-\eqref{J1-estim-part3-unc}, for
$N^{7/12 + \eps} \le H \le \odi{N}$ we finally obtain
\begin{equation}
\label{I2-estim-final-part3-unc}
I_{2} \ll H \exp \Big( - \frac{c_{1}}{4}  \Big( \frac{L}{\log L} \Big)^{1/3} \Big) .
\end{equation}

Now we estimate $I_{3}$. By \eqref{UH-estim}, Lemma \ref{zac-lemma-series} 
and a partial integration argument
we get
\begin{align}
\label{I3-estim-final-part3-unc}
I_3
\ll
\int_{B/H}^{1/2}  
\vert \Stilde_2(\alpha)\vert^{2} \frac{\dx \alpha}{\alpha}
\ll
  N^{1/2}  L +\frac{HL^{2}}{B}+
L
\int_{B/H}^{1/2}  
(\xi N^{1/2}  + L) \frac{\dx \xi}{\xi^2}
\ll
\Bigl( N^{1/2}  +\frac{H}{B}\Bigr) L^{2} .
\end{align}
Now using \eqref{approx-th2-part3}-\eqref{I1-eval-part3-unc}
and \eqref{I2-estim-final-part3-unc}-\eqref{I3-estim-final-part3-unc},
and choosing $0<c<c_{1}/4$ in \eqref{B-def},
we have that there exists a constant $C=C(\eps)>0$ such that

\begin{equation*}
\sum_{n = N+1}^{N + H}  e^{-n/N}
r^{\second}_{2,2}(n)
 =
 \frac{\pi H}{4e} 
 + 
\Odig{ H \exp \Big( - C  \Big( \frac{L}{\log L} \Big)^{1/3} \Big) + \frac{H^{2}}{N}}
\end{equation*}
uniformly for     $N^{7/12+\eps} \le H \le \odi{N}$.
Theorem \ref{asymp-unc-part3}  hence follows for $N^{7/12+\eps} \le H \le N^{1-\eps}$
since the exponential weight $e^{-n/N}$
can be removed as we did at the bottom of the proof of Theorem \ref{asymp-part3}.
\qed 
 
\begin{Remark}
\label{Remark-dopo-teo2}
Using the finite-sum approach we need to define 
$T_{2}(\alpha) = \sum_{1 \le m^{2} \le  N} e(m^{2} \alpha )$
and $f_{2}(\alpha) =(1/2)\sum_{1 \leq m\leq N} m^{-1/2}e(m\alpha)$.
Theorem 4.1 of Vaughan \cite{Vaughan1997} gives
$\vert T_{2}(\alpha) - f_{2}(\alpha) \vert \ll (1+\vert \alpha \vert N)^{1/2}$.
The main term comes from the integral of $f_{2}(\alpha)^{2}U(-\alpha,H)$
but we also need to evaluate the quantity
\[
\Bigl\vert   
\int_{-B/H}^{B/H}  (T_{2}(\alpha)^{2}-f_{2}(\alpha)^{2})U(-\alpha,H)
 e(-N\alpha)\, \dx \alpha
 \Bigr\vert 
 \ll
 \frac{NB^{1/2}}{H^{1/2}}.
\]
Since the expected order of magnitude of the 
main term is $H$, the previous estimate
is under control  if and only if $H\ge N^{2/3} B^{1/3}$ which is weaker 
than the result we obtain.
Similar remarks apply for the other problems studied in the remaining sections.
 \end{Remark}

\section{Proof of Theorem \ref{asymp-part4} }  
  
%
Letting $1<A=A(N)<H/2$ to be chosen later, by \eqref{main-defs} 
and \eqref{omega-def}-\eqref{omega-approx}
it is an easy matter to see that
\begin{align}
\notag
\sum_{n = N+1}^{N + H} & e^{-n/N}
r^{\prime}_{2,2}(n) 
= 
\int_{-1/2}^{1/2}  \Ttilde_2(\alpha) \omega(\alpha)  U(-\alpha,H) e(-N\alpha) \, \dx \alpha
\\
\notag
 &=
  \int_{-1/2}^{1/2}  (\Ttilde_2(\alpha)- \Stilde_2(\alpha)) \omega(\alpha)  U(-\alpha,H) e(-N\alpha) \, \dx \alpha
 \\
 \notag
 & \hskip0.3cm 
 + \int_{-A/H}^{A/H}   \Bigl(   \frac{\pi}{4z} -   \frac{\pi^{1/2}}{4z^{1/2}}\Bigr)
 U(-\alpha,H) e(-N\alpha)\, \dx \alpha
  +
 \int_{-A/H}^{A/H} 
 \Etilde_{2}(\alpha) 
 \omega(\alpha) U(-\alpha,H)
 e(-N\alpha)\, \dx \alpha 
\\
 \notag
 & \hskip0.3cm
 +
\int_{-A/H}^{A/H}   \frac{\pi}{2z}
 \Bigl( \sum_{\ell=1}^{+\infty}   e^{-\ell^{2}\pi^{2}/z} \Bigr)
 U(-\alpha,H) e(-N\alpha)\, \dx \alpha 
+ 
\int\limits_{\mathclap{[-1/2,-A/H]\cup  [A/H,1/2]}} 
\Stilde_2(\alpha)\omega(\alpha)   U(-\alpha,H)  e(-N\alpha) \, \dx \alpha
\\
\label{approx-th1-part4} 
&
= I_{0} + I_{1}+I_{2}+I_{3} +I_{4}, 
\end{align}
say.
Using Lemma \ref{trivial-lemma}  and the  Cauchy-Schwarz inequality we have
\begin{align*} 
I_{0} 
\ll
N^{1/4}
\Bigl(
\int_{-1/2}^{1/2}
      \vert \omega(\alpha)  \vert^{2}
        \vert U(\alpha,H)\vert 
     \, \dx \alpha
     \Bigr)^{1/2}
 \Bigl(
      \int_{-1/2}^{1/2}
       \vert U(\alpha,H) \vert
        \, \dx \alpha        
\Bigr)^{1/2}.
\end{align*}
By Lemma \ref{zac-lemma-series}, \eqref{UH-estim} 
and a partial integration argument we obtain
\begin{equation}
\label{I0-estim-omega}
I_{0}
\ll
N^{1/4}
(N^{1/2}L  + H L)^{1/2}L^{1/2}  
\ll
N^{1/2}L+ H^{1/2}N^{1/4}L.
\end{equation}

Now we evaluate $I_{1}$. Using Lemma~4 of \cite{LanguascoZ2016a} 
and \eqref{UH-estim}
we immediately get 
\begin{equation}
\label{I1-eval-part4}
I_{1}
=
 \sum_{n = N+1}^{N + H}  \Bigl(\frac{\pi}{4}   - \frac{1}{4n^{1/2}}   \Bigr)  e^{-n/N} + \Odig{\frac{H}{N}}
+
\Odig{\int_{A/H}^{1/2} \frac{\dx\alpha}{\alpha^2}}
  =
\frac{\pi H }{4e} + \Odig{\frac{H}{N^{1/2}}+\frac{H^2}{N}+\frac{H}{A}}.
\end{equation}
To have that $\pi H/(4e)$ dominates in $I_{0} + I_{1}$ we need that 
$A\to \infty$, $H=\odi{N}$ and  $H=\infty(N^{1/2}L^{2})$.

Now we estimate $I_{3}$. Assuming $H=\infty(N^{1/2}A)$, 
by \eqref{omega-Y-estim}-\eqref{z-estim}, we have
\begin{align}
\notag
I_{3}
&\ll
\frac{HN}{e^{\pi^2 N}} \int_{-1/N}^{1/N}  \dx \alpha
+
\frac{H}{e^{\pi^2 H^2/N}} \int_{1/N}^{1/H}  \frac{\dx \alpha}{\alpha} 
+
\int_{1/H}^{A/H}  \frac{\dx \alpha}{\alpha^2e^{\pi^2/(N\alpha^2)}}  
\\
\label{I3-estim-final-part4}
&
\ll
\frac{H}{e^{\pi^2 N}} + \frac{HL}{e^{\pi^2 H^2/N}} + \frac{H}{e^{\pi^2 H^2/(NA^{2})}}
\end{align}
which is $\odi{H}$
provided that  $H=\infty(N^{1/2} \log L)$ and $H=\infty(N^{1/2}  A )$.

Now we estimate $I_{2}$. Recalling  $H=\infty(N^{1/2}A)$,
for every $\vert \alpha \vert \le A/H$
we have, by \eqref{omega-approx}-\eqref{omega-Y-estim}, 
that $\vert\omega(\alpha) \vert \ll \vert z \vert ^{-1/2}$. 
Hence
\[
I_{2} 
\ll
 \int_{-A/H}^{A/H}
      \vert \Etilde_{2}(\alpha)  \vert  \frac{\vert U(\alpha,H)\vert}{\vert z \vert ^{1/2}}
     \, \dx \alpha.
\]
Using  \eqref{z-estim} and the Cauchy-Schwarz inequality and Lemma~3
of \cite{LanguascoZ2016a} we get 
\begin{align}
\notag
I_{2}
&\ll 
H N^{1/2} \Bigl(  \int_{-1/N}^{1/N}\!\!\!\! \dx \alpha  \Bigr)^{1/2}
\Bigl(
 \int_{-1/N}^{1/N} \!\! \vert \Etilde_{2}(\alpha)  \vert^{2}   \, \dx \alpha 
 \Bigr)^{1/2}
\!\!\! +
 H \Bigl(  \int_{1/N}^{1/H} \!\!  \frac{\dx \alpha}{\alpha^{1/2}}\Bigr)^{1/2} 
\Bigl(
 \int_{1/N}^{1/H}\!\!  \vert \Etilde_{2}(\alpha)  \vert^{2}  
 \frac{\dx \alpha}{\alpha^{1/2}}
  \Bigr)^{1/2}
  \\
  \notag
     &\hskip2cm 
     +
  \Bigl(  \int_{1/H}^{A/H}   \frac{\dx \alpha}{\alpha^{3/2}}\Bigr)^{1/2} 
\Bigl(
 \int_{1/H}^{A/H}  \vert \Etilde_{2}(\alpha)  \vert^{2}   
 \frac{\dx \alpha}{\alpha^{3/2}}
   \Bigr)^{1/2}
 \\
 \notag 
 &\ll 
H N^{-1/4} L 
+ 
H^{3/4}  N^{1/4}  L  \Bigl(  \frac{1}{H}+ \int_{1/N}^{1/H}  \frac{\dx \xi}{\xi^{1/2}}   \Bigr)^{1/2} 
 +
H^{1/4} N^{1/4} L \Bigl( H^{1/2} +   \int_{1/H}^{A/H}  \frac{\dx \xi}{\xi^{3/2}}   \Bigr)^{1/2}
 \\
 \label{I2-estim-final-part4}
 &\ll  
 H^{1/2}N^{1/4} L.
\end{align}

Remark that $I_2=\odi{H}$ provided that $H=\infty(N^{1/2} L^{2})$.

Now we estimate $I_{4}$. By \eqref{UH-estim}, Lemma \ref{zac-lemma-series} 
and a partial integration argument
we get
\begin{align}
\notag
I_4
&\ll
\int_{A/H}^{1/2}
\vert  \Stilde_2(\alpha) \omega(\alpha) \vert
\frac{\dx \alpha}{\alpha}
\ll
\Bigl(
\int_{A/H}^{1/2}  
\vert \Stilde_2(\alpha)\vert^2 \frac{\dx \alpha}{\alpha}
\Bigr)^{1/2}
\Bigl(
\int_{A/H}^{1/2}  
\vert \omega(\alpha) \vert^2\frac{\dx \alpha}{\alpha}
\Bigr)^{1/2}
\\
\notag
&\ll
\Bigl(
 N^{1/2}   L
+
\frac{HL^{2}}{A}
+
L
\int_{A/H}^{1/2}  
(\xi N^{1/2}  + L) \frac{\dx \xi}{\xi^2}
\Bigr)^{1/2}
\Bigl(
 N^{1/2} 
+
\frac{HL}{A}
+
\int_{A/H}^{1/2}  
(\xi N^{1/2}  + L) \frac{\dx \xi}{\xi^2}
\Bigr)^{1/2}
\\
\label{I4-estim-final-part4}
&\ll
 L^{3/2} \Bigl(
 N^{1/2}
+
\frac{H}{A}
\Bigr) 
\end{align}
which is $\odi{H}$
provided that $A=\infty(L^{3/2})$ and $H=\infty(N^{1/2}L^{3/2})$.

Combining the conditions on $H$ and $A$ we can choose $A=L^{2}/(\log L)$
and $H=\infty (N^{1/2}L^{2})$. Hence using \eqref{approx-th1-part4}-\eqref{I4-estim-final-part4} we  
can write  
\begin{equation*}
\sum_{n = N+1}^{N + H}  e^{-n/N}
r^{\prime}_{2,2}(n)   
 =
\frac{\pi H }{4e}  + \Odig{\frac{H^2}{N}+ \frac{H \log L} {L^{1/2}} +  N^{1/2} L^{3/2} +
H^{1/2}N^{1/4} L}.
\end{equation*}
Theorem \ref{asymp-part4}  follows for $\infty(N^{1/2}L^{2}) \le H \le \odi{N}$
since the exponential weight $e^{-n/N}$
can be removed as we did at the bottom of the proof of Theorem \ref{asymp-part3}.
\qed

\section{Proof of Theorem \ref{asymp-unc-part4}}  

By  \eqref{main-defs} 
and \eqref{omega-def}-\eqref{omega-approx}, it is an easy matter to see that
\begin{align}
\notag
\sum_{n = N+1}^{N + H} &  e^{-n/N}
r^{\prime}_{2,2}(n) 
= 
\int_{-1/2}^{1/2}  \Ttilde_2(\alpha) \omega(\alpha)  U(-\alpha,H) e(-N\alpha) \, \dx \alpha
\\
\notag
 &=
  \int_{-1/2}^{1/2}  (\Ttilde_2(\alpha)- \Stilde_2(\alpha)) \omega(\alpha)  U(-\alpha,H) e(-N\alpha) \, \dx \alpha
 \\
\notag
&\hskip0.3cm
+
 \int_{-B/H}^{B/H}      \Bigl(   \frac{\pi}{4z} -   \frac{\pi^{1/2}}{4z^{1/2}}\Bigr) 
  U(-\alpha,H) 
  e(-N\alpha)\, \dx \alpha 
  +
 \int_{-B/H}^{B/H}
\Etilde_{2}(\alpha) 
 \omega(\alpha)   U(-\alpha,H) 
 e(-N\alpha)\, \dx \alpha 
\\ 
\notag
 & \hskip0.3cm
  +
 \int_{-B/H}^{B/H}
  \frac{\pi}{2z}
 \Bigl( \sum_{\ell=1}^{+\infty}   e^{-\ell^{2}\pi^{2}/z} \Bigr)
 U(-\alpha,H) 
 e(-N\alpha)\, \dx \alpha 
 +
\int\limits_{\mathclap{[-1/2,-B/H]\cup  [B/H,1/2]}} 
\Stilde_2(\alpha)\omega(\alpha)   U(-\alpha,H)  e(-N\alpha) \, \dx \alpha
\\
\label{approx-th2-part4} 
&
= I_{0}+ I_{1}+I_{2}+I_{3}+I_{4}   , 
\end{align}
say, where $B$ is defined in \eqref{B-def}. 
$I_{0}$ can be estimated as in \eqref{I0-estim-omega} and gives
\begin{equation}
\label{I0-estim-omega-unc}
I_{0}
\ll
N^{1/2}L+ H^{1/2}N^{1/4}L.
\end{equation}

$I_{1}$ can be evaluated as in \eqref{I1-eval-part4} and we get
\begin{equation}
\label{I1-eval-part4-unc}
I_{1}
   =
\frac{\pi H}{4e} + \Odig{\frac{H^{2}}{N}+\frac{H}{B}}.
\end{equation}

Now we estimate $I_{2}$.
Using  \eqref{UH-estim} and  the Cauchy-Schwarz inequality we obtain 
\begin{align}
\label{I2-estim1-part4-unc}
I_{2}
\ll
H
\Bigl( \int_{-B/H}^{B/H} 
      \vert \Etilde_{2}(\alpha) \vert^{2}    
     \, \dx \alpha
\Bigr)^{1/2}
\Bigl(
 \int_{-B/H}^{B/H} 
       \vert \omega(\alpha)  \vert^{2}  
        \, \dx \alpha 
\Bigr)^{1/2} 
    =
    H   (J_{1}J_{2})^{1/2},
\end{align}
say.
Using  Lemma \ref{LP-Lemma-gen-unc}   we can write
\begin{align}
 \label{J1-estim-part4-unc}
J_{1} 
\ll
\exp \Big( - c_{1}  \Big( \frac{L}{\log L} \Big)^{1/3} \Big) 
\end{align}
provided that
$\infty(1/N)<B/H < N^{-7/12 - \eps/2}$, i.e. 
$N^{7/12 + \eps} \le H \le \odi{N}$ suffices.

Using  Lemma \ref{zac-lemma-series} with $\ell=2$ we have
\begin{equation}
\label{J2-estim-part4-unc}
 J_{2}
 \ll  
 \frac{N^{1/2} B}{H} + L \ll L.
\end{equation}

Combining \eqref{I2-estim1-part4-unc}-\eqref{J2-estim-part4-unc}
for $N^{7/12 + \eps} \le H \le \odi{N}$  we finally obtain
\begin{equation}
\label{I2-estim-final-part4-unc}
I_{2} \ll H \exp \Big( - \frac{c_{1}}{4}  \Big( \frac{L}{\log L} \Big)^{1/3} \Big) .
\end{equation}

Now we estimate $I_{3}$.
By  \eqref{omega-Y-estim}-\eqref{z-estim}, we have
\begin{equation}
\label{I3-estim-final-part4-unc}
I_{3}
\ll
\frac{HN}{e^{\pi^2 N}} \int_{-1/N}^{1/N} \dx \alpha
+
\frac{H}{e^{\pi^2 H^2/(NB^{2})}} \int_{1/N}^{B/H}  \frac{\dx \alpha}{\alpha}  
\ll
 H \exp \Big( - \frac{c_{1}}{4}  \Big( \frac{L}{\log L} \Big)^{1/3} \Big),
\end{equation}
since $N^{7/12 + \eps} \le H \le \odi{N}$.

$I_{4}$ can be estimated as in \eqref{I4-estim-final-part4} and gives
\begin{equation}
\label{I4-estim-final-part4-unc}
I_4
\ll
 L^{3/2} \Bigl(
 N^{1/2}
+
\frac{H}{B}
\Bigr) .
\end{equation}

Now using \eqref{approx-th2-part4}-\eqref{I1-eval-part4-unc}
and \eqref{I2-estim-final-part4-unc}-\eqref{I4-estim-final-part4-unc},
and choosing $0<c<c_{1}$ in \eqref{B-def},
we have that there exists a constant $C=C(\eps)>0$ such that
\begin{equation*} 
\sum_{n = N+1}^{N + H}  e^{-n/N}
r^{\prime}_{2,2}(n)
 =
 \frac{\pi H}{4e} 
 + 
\Odig{  H \exp \Big( - C  \Big( \frac{L}{\log L} \Big)^{1/3} \Big) +\frac{H^{2}}{N}}
\end{equation*}
uniformly for 
for    $N^{7/12+\eps} \le H \le \odi{N}$. 
Theorem \ref{asymp-unc-part4}  hence follows for $N^{7/12+\eps} \le H \le N^{1-\eps}$ 
since the exponential weight $e^{-n/N}$
can be removed as we did at the bottom of the proof of Theorem \ref{asymp-part3}.
\qed

\renewcommand{\bibliofont}{\normalsize}

\bigskip
\noindent
Alessandro Languasco, Dipartimento di Matematica, Universit\`a
di Padova, Via Trieste 63, 35121 Padova, Italy. {\it e-mail}: languasco@math.unipd.it

\medskip
\noindent
Alessandro Zaccagnini, Dipartimento di Matematica e Informatica, Universit\`a di Parma, Parco
Area delle Scienze 53/a, 43124 Parma, Italy. {\it e-mail}: alessandro.zaccagnini@unipr.it
\end{document}